\documentclass{amsart}
\usepackage{graphicx}

\twocolumn \textwidth=16.5cm\textheight=24.6cm
\begin{document}

\noindent
{\LARGE\bf A Mathematical Theory for}

\smallskip
\noindent
{\LARGE\bf Random Solid Packings}

\vspace{0.5cm}
\noindent{\bf Chuanming Zong}\footnote{Email address: cmzong@math.pku.edu.cn}

\medskip
\noindent
School of Mathematical Sciences, Peking University, Beijing
100871, China.

\vspace{0.8cm}
\noindent
{\bf Packings of identical objects have fascinated both scientists and laymen alike for centuries, in particular the sphere packings and the packings of identical regular tetrahedra. Mathematicians have tried for centuries to determine the densest packings; Crystallographers and chemists have been fascinated by the lattice packings for centuries as well. On the other hand, physicists, geologists, material scientists and engineers have been challenged by the mysterious random packings for decades. Experiments have shown the existence of the dense random sphere packings and the loose random sphere packings for more than half a century. However, a rigorous definition for them is still missing. The purpose of this paper is to review the random solid packings and to create a mathematical theory for it.}

\bigskip
\noindent
dense random packing $|$ loose random packing $|$ Platonic solids

\vspace{0.6cm}
\noindent
More than 2300 years ago, Aristotle (384--322 BC) stated: \lq\lq {\it Among surfaces it is agreed that there are three figures which fill the place that contain them -- the triangle, the square and the hexagon;
among solids only two, the pyramid and the cube.}" [1, 306b, p. 319]. Here \lq\lq pyramid" means the regular tetrahedron, a Platonic solid. Thus, in modern terminology, Aristotle's assertion can be taken to mean: {\it The space can be tiled by congruent regular tetrahedra.}

In fact, Aristotle was wrong. Identical regular tetrahedra cannot tile the whole space!
According to Struik \cite{stru25}, the incorrectness of Aristotle's statement was likely known to Johannes M\"uller (1436-1476)
and a definite observation of the incorrectness was made by Paulus van Middelburg (1445-1534). Then a natural problem arose: {\it What is the maximum density, when one pack small identical regular tetrahedra in a big box}?

1594, Sir Walter Raleigh (1552--1618) asked his assistant Thomas Harriot (1560--1621) to figure out the best way to pack identical cannonballs. As a result, he discovered the {\it face-centered cubic packing} with density $\pi /\sqrt{18}$. Harriot was an ardent atomist. He believed that the secrets of the universe were to be revealed through the patterns and packings of small spherical atoms.

Based on Harriot's discovery, in 1611 Kepler \cite{kepl11} made the following conjecture:
{\it When one pack small identical spheres in a big container, the maximal density is $\pi /\sqrt{18}$.}

\medskip
In his talk presented at the ICM 1900 in Paris, Hilbert \cite{hilb01} proposed 23 unsolved mathematical problems. As a part of his 18th problem, based on Aristotle's mistake and Kepler's conjecture, he wrote: \lq\lq {\it I point out the following question, related to the preceding one, and important to number theory and perhaps sometimes useful to physics and chemistry: How can one arrange most densely in space an infinite number of equal solids of given form, e.g., spheres with given radii or regular tetrahedra with given edges $($or in prescribed position$)$, that is, how can one so fit them together that the ratio of the filled to the unfilled space may be as great as possible}?"

In the history, packing problems have been studied by many great mathematicians, including Newton, Gauss, Hermite, Fedorov, Voronoi and Minkowski. The spherical case was solved by Hales \cite{hale05} only in 2005, with the assistance of a complicated computer programme; The tetrahedral case is still widely open.  For mathematical achievements on sphere packings and tetrahedron packings we refer to [5-11] and their references. The purpose of this paper is to review random solid packings, the approaches made by physicists, material scientists and engineers, and to create a mathematical theory for it.

\vspace{0.4cm}

\noindent {\large\bf Random Ball Packings}

\smallskip\noindent
Being suspicious of Kepler's conjecture, in 1958 Coxeter \cite{coxe58} wrote:\lq\lq {\it The densest lattice-packing is not necessarily the densest packing.}" Therefore, random packing is a reasonable attempt to uncover a counterexample, if there were such examples. On the other hand, in 1959 Bernal \cite{bern59} suggested to use geometrical models to study the structure of liquids. In particular, a random packing of equal balls may provide a useful model for an ideally simple liquid. For these motives, in 1960 Scott \cite{scot60} made a simple but extremely important experiment.

When a lot of identical small balls were poured into a graduated cylinder, shaking and tapping quickly reduced the volume. Continued gentle shaking produced no further reduction in volume. More vigorous shaking increased the volume. On the other hand, when the container was tipped horizontally and then slowly rotated around its axis and gradually returned to the vertical position the volume increased to a certain level. For all the rigid containers used, a similar effect was observed, namely, there is a range of random packing densities lying between two well-defined limits. The limits are called {\it dense random packing density $\delta_d(B)$} and {\it loose random packing density $\delta_l(B)$}, respectively.

Scott used small identical steel balls and two kinds of containers, one spherical and the other cylindrical, for his experiment. He obtained very stable results: {\it in both cases, the dense random packing density is about $0.63$ and the loose random packing density is about $0.59.$}

It is quite surprising that the densities of the dense random sphere packings (as well as the loose random packings) are very stable; Even more surprising is that the stable density is rather far away from $\pi /\sqrt{18}\approx 0.7405\ldots ,$ the density of the densest sphere packings (see \cite{hale05}).

\smallskip
Suspected that the specific gravity of the balls may have had an influence on the arrangement of the balls in random packings, and therefore the packing densities, in 1962 Rutgers \cite{rutg62} made an experiment on random packing of nylon balls. Indeed, he obtained stable dense random packings with smaller densities. Furthermore, he also observed slight effects by the filling time and filling height as shown in Table 1.

\bigskip
\centerline{
\begin{tabular}{|c|c|c|c|}
\hline $\begin{array}{c}
{\rm filling} \\
{\rm time\ (sec.) }
\end{array}$& $\begin{array}{c}
{\rm filling} \\
{\rm height\ (cm.) }
\end{array}$& $\delta_l(B)$ & $\delta_d(B)$ \\
\hline $1.5$ & $17.5$ & $0.556$ & $0.588$\\
\hline $1.5$ & $8.7$ & $0.568$ & $0.597$ \\
\hline $11$ & $17.5$ & $0.591$ & $0.598$\\
\hline $11$ & $8.7$ & $0.583$ & $0.601$ \\
\hline
\end{tabular}}

\medskip
\noindent{\footnotesize {\bf Table 1.} Random nylon ball packings with different conditions.}

\medskip
Further similar experiments on random ball packings were reported by several authors [16-18]. In particular, Dong and Ye reported random ball packings with different materials, as listed in Table 2.

\bigskip
\noindent
\begin{tabular}{|c|c|c|c|c|c|}
\hline material & lead & steel & glass & plastic & plastic \\
\hline $g/cm^3$ & $11.30$ & $7.85$ & $2.5$ & $1.032$ & $0.912$ \\
\hline $\hspace{-0.2cm}
\begin{array}{c}
\delta_d(B)\\
{\rm errors}\ (\pm )\\
{\rm exp.\ times}
\end{array}\hspace{-0.2cm}$&
$\hspace{-0.2cm}
\begin{array}{c}
0.641\\
\cdot\\
\cdot
\end{array}\hspace{-0.2cm}$&
$\hspace{-0.2cm}
\begin{array}{c}
0.639\\
0.003\\
24
\end{array}\hspace{-0.2cm}$&
$\hspace{-0.2cm}
\begin{array}{c}
0.600\\
0.002\\
10
\end{array}\hspace{-0.2cm}$ &
$\hspace{-0.2cm}
\begin{array}{c}
0.641\\
0.002\\
12
\end{array}\hspace{-0.2cm}$ &
$\hspace{-0.2cm}
\begin{array}{c}
0.640\\
0.001\\
10
\end{array}\hspace{-0.2cm}$ \\
\hline $\hspace{-0.2cm}
\begin{array}{c}
\delta_l(B)\\
{\rm errors}\ (\pm )\\
{\rm exp.\ times}
\end{array}\hspace{-0.2cm}$&
$\hspace{-0.2cm}
\begin{array}{c}
0.553\\
\cdot\\
\cdot
\end{array}\hspace{-0.2cm}$&
$\hspace{-0.2cm}
\begin{array}{c}
0.553\\
0.004\\
14
\end{array}\hspace{-0.2cm}$&
$\hspace{-0.2cm}
\begin{array}{c}
0.506\\
0.003\\
12
\end{array}\hspace{-0.2cm}$ &
$\hspace{-0.2cm}
\begin{array}{c}
0.581\\
0.002\\
10
\end{array}\hspace{-0.2cm}$ &
$\hspace{-0.2cm}
\begin{array}{c}
0.550\\
0.001\\
10
\end{array}\hspace{-0.2cm}$ \\
\hline
\end{tabular}

\medskip
\noindent{\footnotesize {\bf Table 2.} Random ball packings with different materials.}

\medskip
It is very surprising that $\delta_d(B)$ is not increasing with respect to the specific gravity of the balls. In particular, it is hard to imagine that the density of the dense random packing of glass balls is so much small, comparing to the other cases.

\vspace{0.4cm}

\noindent
{\large\bf Statistical Explanations}

\smallskip\noindent
To study the structure of a random ball packing, Bernal and Mason \cite{bern60} and Scott \cite{scot62} designed some very practical experiments. For example, Scott did the following one. A big number of small steel balls were poured into a cylindrical container. The balls were shaken down to achieve the dense random packing and then after the container was warmed to about $65^\circ C$., molten paraffin wax was poured in and the whole allowed to cool. Take the balls with wax out of the container and remove the outside layers of the balls and wax. Under the comparator the coordinates of the balls in a roughly spherical cluster were determined. This experiment releases two kinds of information, which are important to understand the random packing.

Let $f(x)$ denote the average number of the neighbours within a radial distance of $(x+1)d$ to a ball, where $d$ is the diameter of the balls. In other words, let ${\bf x}_1$, ${\bf x}_2$, $\ldots $, ${\bf x}_m$ denote the centers of the balls in the cluster of the random packing, let $n_i$ denote the number of the points ${\bf x}_j$, ${\bf x}_j\not= {\bf x}_i$, such that $\| {\bf x}_i, {\bf x}_j\|\le (1+x)d$, then
$$f(x)={1\over m}\sum_{i=1}^mn_i.$$
As usual, $\| {\bf x}, {\bf y}\|$ denotes the Euclidean distance between ${\bf x}$ and ${\bf y}$.
Based on Scott's experiment, in 1968 Mason \cite{maso68} obtained a diagram for $f(x)$. The following diagram was obtained by Gotoh and Finney \cite{goto74} based on a 7934 center close-packed model of Finney. It is particular interesting to note that the average touching number of the balls in the dense random packing is less than six, which is much smaller than twelves, the maximal number of identical balls which can be brought into contact with a fixed one (see \cite{zong99}).

\begin{figure}[ht]
\centering
\includegraphics[height=4.5cm,width=5.5cm,angle=0]{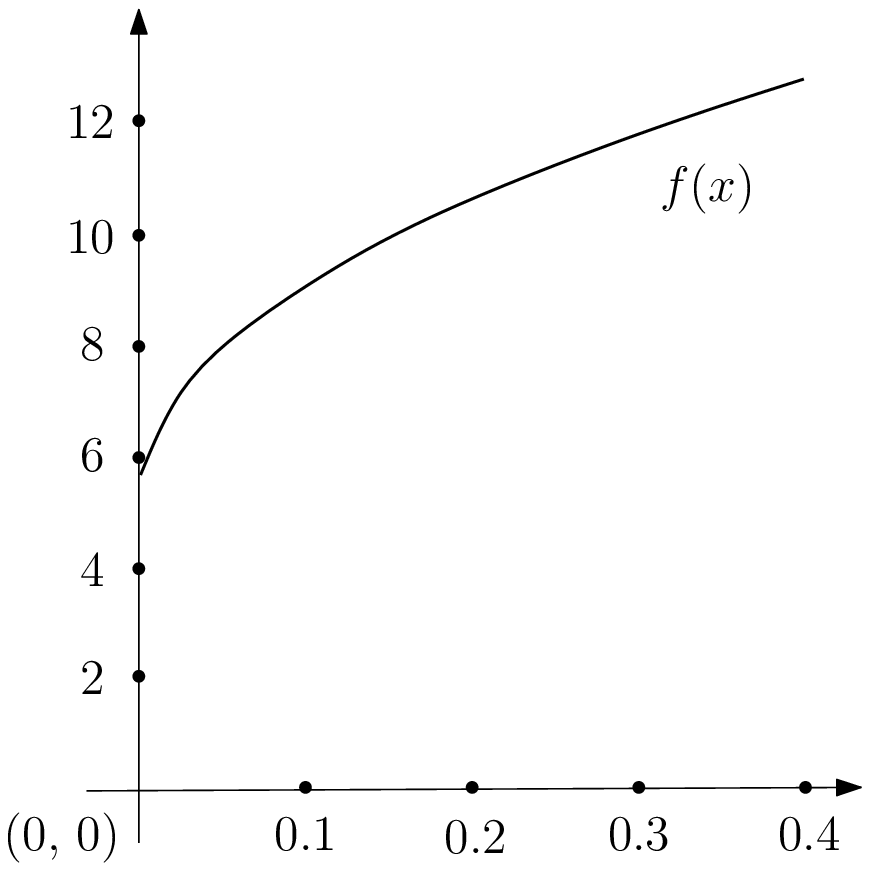}
\end{figure}

\smallskip
\noindent{\footnotesize {\bf Fig. 1.} Average radial distance distributions in dense random ball packings.}

\medskip
In the random packing, every ball is associated with a {\it Voronoi cell}. When the ball belongs to the cluster, its Voronoi cell is a polytope. According to Finney \cite{finn70}, most of the Voronoi cells in the dense random ball packings have $12\sim 17$ faces (see Table 3), and their average is about $14.251.$

\bigskip
\centerline{
\begin{tabular}{|c|c|c|c|c|c|c|}
\hline
$\begin{array}{c}
{\rm face}\\
{\rm number}
\end{array}$ & 12 & 13 & 14 & 15 & 16 & 17\\
\hline fraction (\%) & $4.2$ & $20.5$ & $35.2$ & $27.3$ & $10.5$ & $1.6$\\
\hline
\end{tabular}}

\medskip
\noindent{\footnotesize {\bf Table 3.} Distribution of the face numbers of the Voronoi cells in the dense random ball packings.}

\medskip
Based on these facts, Gotoh and Finney \cite{goto74} studied the most probable tetrahedron in the dense random ball packings, with ball centers as its vertices and three unit edges (length $d$) with a common vertex. By computing the fraction contained by the four balls at the vertices, they were able to achieve an explanation for the density $\delta_d(B)\approx 0.6357$. There are other mathematical verifications and explanations for this fact (see [24-27]). Nevertheless, a mathematical proof is still missing.

\vspace{0.4cm}
\noindent {\large\bf Random Packings of Platonic Solids}

\smallskip\noindent
{\it If the balls in the random packings are replaced by other identical objects, such as the Platonic solids, what will happen}?

In 1993, Dong and Ye \cite{dong93} studied the cases for cubes and regular tetrahedra. It is known to everybody that identical cubes can tile the whole space. In other words, the density of the densest cube packings is $1$. However,
with small identical wood cubes to make dense random packings for eight independent times, they obtained a stable density $0.640$ with an error of $\pm 0.003$. For the densities $\delta_d(T)$ of the dense random packings of regular tetrahedra, they tried with different specific gravities and obtained the results listed in Table 4.

\bigskip
\centerline{
\begin{tabular}{|c|c|c|c|c|c|}
\hline $\hspace{-0.2cm}
\begin{array}{c}
{\rm gravities}\\
(g/cm^3)
\end{array}\hspace{-0.2cm}$ &3.78& 2.625& 1.68&0.945 & 0.659 \\
\hline $\hspace{-0.2cm}
\begin{array}{c}
\delta_d(T)\\
{\rm errors}\ (\pm )\\
{\rm exp.\ times}
\end{array}\hspace{-0.2cm}$&
$\hspace{-0.2cm}
\begin{array}{c}
0.494\\
0.004\\
37
\end{array}\hspace{-0.2cm}$&
$\hspace{-0.2cm}
\begin{array}{c}
0.478\\
0.009\\
35
\end{array}\hspace{-0.2cm}$&
$\hspace{-0.2cm}
\begin{array}{c}
0.461\\
0.006\\
38
\end{array}\hspace{-0.2cm}$ &
$\hspace{-0.2cm}
\begin{array}{c}
0.500\\
0.001\\
40
\end{array}\hspace{-0.2cm}$ &
$\hspace{-0.2cm}
\begin{array}{c}
0.489\\
0.006\\
32
\end{array}\hspace{-0.2cm}$ \\
\hline $\hspace{-0.2cm}
\begin{array}{c}
\delta_l(T)\\
{\rm errors}\ (\pm )\\
{\rm exp.\ times}
\end{array}\hspace{-0.2cm}$&
$\hspace{-0.2cm}
\begin{array}{c}
0.387\\
0.003\\
38
\end{array}\hspace{-0.2cm}$&
$\hspace{-0.2cm}
\begin{array}{c}
0.388\\
0.001\\
24
\end{array}\hspace{-0.2cm}$&
$\hspace{-0.2cm}
\begin{array}{c}
0.385\\
0.007\\
14
\end{array}\hspace{-0.2cm}$ &
$\hspace{-0.2cm}
\begin{array}{c}
0.372\\
0.005\\
24
\end{array}\hspace{-0.2cm}$ &
$\hspace{-0.2cm}
\begin{array}{c}
0.364\\
0.005\\
33
\end{array}\hspace{-0.2cm}$ \\
\hline
\end{tabular}}

\medskip
\noindent{\footnotesize {\bf Table 4.} Random packings with regular tetrhedra of different specific gravities.}

\medskip
Similar to the ball case, the density $\delta_d(T)$ of the dense random tetrahedron packings is not monotonic
with respect to the specific gravity.

\medskip
In 2010, Jaoshvili, Esakia, Porrati and Chaikin \cite{jaos10} reported a density $0.76\pm 0.02$ of random
tetrahedral dice packings. However, according to Baker and Kudrolli \cite{bake10}, the protocol by which the packings were prepared was not clear. In fact, Baker and Kudrolli studied the dense random packings of all Platonic solids (as shown in Table 5). The actual particles studied in their experiment are plastic dices (specific gravity of $1.16$ $g/cm^3$) which have slightly rounded edges.

\bigskip
\centerline{
\begin{tabular}{|c|c|c|}
\hline shape ($P$) &  $\delta_d(P)$ & $\delta_l(P)$\\
\hline tetrahedral dice &  $0.64\pm 0.01$ & $0.54\pm 0.01$\\
\hline cubic dice &  $0.67\pm 0.02$ & $0.57\pm 0.01$\\
\hline octahedral dice & $0.64\pm 0.01$ & $0.57\pm 0.01$\\
\hline dodecahedral dice & $0.63\pm 0.01$ & $0.56\pm 0.01$\\
\hline icosahedral dice & $0.59\pm 0.01$ & $0.53\pm 0.01$\\
\hline
\end{tabular}}

\medskip
\noindent{\footnotesize {\bf Table 5.} Random packings with Platonic solid dices.}

\medskip
It is interesting to note that, although the icosahedral dice is the closest to the ball in shape among all the five Platonic solids, their density difference is the biggest.

\vspace{0.4cm}

\noindent {\large\bf A Mathematical Theory for Random Solid Packings}

\smallskip\noindent
In 1960, Bernal and Mason \cite{bern60} wrote: \lq\lq {\it The figure for the occupied volume of random close packing - $0.64$ - must be mathematically determinable, although so far as we know undetermined.}" More than half a century later, a mathematically rigorous definition of the random close packing still remains elusive (see [26, 30]).

\medskip
Random packings is a fundamental problem in nature, which attracts the attention of both scientists and engineers. Therefore, it deserves a rigorous definition and a systematic study. Since it is a phenomenon under certain physical conditions such as shaking and tapping, the expected definition cannot be purely mathematical.

As usual, we call a set $Q$ convex if it contains all the segments with ends belonging to $Q$. Let $K$ denote a convex body in our space $\mathbb{E}^3$, a bounded convex set with boundary. Clearly, all balls and the Platonic solids are convex bodies. If we make a hole in a ball, it is no longer convex.

Let $\Im $ denote a set of rigorously stated experiment conditions and processes to produce the dense random packings. For example, one can set $\Im $ to be \lq\lq pour the objects from $\alpha $ cm above the experiment container in a speed of $\beta $ pieces per second, when the container is full put it on the shaking and tapping machine for $\theta$ hour." Of course, the machine is well-designed for the experiment purpose. Similarly, let $\wp $ denote a set of rigorously stated experiment conditions and processes to produce the loose random packings.

Assume that we have a lot of identical solids with shape $K$ and specific gravity $g$. Let $C_r$ denote a big cylinder with radius $r$ and height $2r$. Fill the cylinder with the solids under the instructions of $\Im$, shake and tap the cylinder following the processes of $\Im$, then by measuring the filled level of the cylinder and the total volume of the solids in the container one obtains a density $\delta^1( K, g, \Im, r)$. Repeat this experiment, we get a sequence of densities
$$\delta^i(K, g, \Im , r), \qquad i=1, 2, 3, \cdots .$$
In fact, we can define $\delta^i(K, g, \Im , r)$ to be the densities of packing identical solids with shape ${1\over r}K$ and specific gravity $g$ in the fixed container $C_1$.

This sequence should be very stable. However, nobody can guarantee it. Let $E(K, g, \Im, r)$ denote the density expectation of this experiment. In other words,
$$E(K, g, \Im , r)=\lim_{n\to \infty}{1\over n}\sum_{i=1}^n\delta^i(K, g, \Im , r).$$
Then, the density of the dense random packings (under the conditions of $\Im $) of $K$ with gravity $g$ is defined by
$$\delta_d(K, g, \Im )=\lim_{r\to \infty}E(K, g, \Im , r).$$

Similarly, one can define $\delta^i(K, g, \wp , r)$, $E(K, g, \wp, r)$ and the density
$\delta_l(K, g, \wp )$ of the loose random packings of $K$ with specific gravity $g$.

The restriction on the packing solids to be convex is not absolutely necessary. However, it simplifies the situation.

Next, we discuss $\delta_d(K, g, \Im )$ and $\delta_l(K,g, \wp)$ from mathematical point of views.

\medskip
If the dense random packing phenomenon indeed does happen for $K$ under the conditions of $\Im $, then the density $\delta_d(K, g, \Im )$ is well-defined. It would be particular interesting to determine the exact values for the ball, the cube, and the regular tetrahedron. On the other hand, if the limit $\delta_d(K, g, \Im )$ does not exist, the so called dense random packings of $K$ will not stable. In this case, to study
$$\overline{\delta_d}(K, g, \Im )=\limsup_{r\to \infty}E(K, g, \Im , r)$$
will turn to be important and interesting.

Similarly, if the loose random packing phenomenon indeed does happen for $K$ under the conditions of $\wp $, then the density $\delta_l(K, g, \wp)$ is well-defined. It would be particular interesting to determine its values for the ball, the cube, and the regular tetrahedron. On the other hand, if the limit $\delta_l(K, g, \wp)$ does not exist, the so called loose random packings of $K$ will not stable. In this case, to study
$$\underline{\delta_l}(K, g, \wp)=\liminf_{r\to \infty}E(K, g, \wp , r)$$
will turn to be interesting.

\medskip
In the definition of the density $\delta (K)$ of the densest packing of $K$, cubes are used as containers (see \cite{zong99}). By approximation, one can easily show that the cubic containers might be replaced by any fixed convex shape. In other words, $\delta (K)$ is independent of the shapes of the containers.

In the definition of the random packing density $\delta_d(K, g, \Im )$, cylinders $C_r$ are used as containers. Random packings were observed and studied as natural phenomena. Therefore it should have certain generality. Namely, the boundaries of the containers should not have influence on $\delta_d(K, g, \Im )$. Let $\widetilde{C_r}$ to be any container satisfying $C_r\subseteq \widetilde{C_r} \subseteq C_{r+d}$, where $d$ is the diameter of $K$. Replacing the containers $C_r$ by $\widetilde{C_r}$ in the definitions of $\delta^i(K, g, \Im , r)$, $E(K, g, \Im, r)$ and $\delta_d(K, g, \Im )$, we obtain corresponding $\widetilde{\delta^i}(K, g, \Im , r)$, $\widetilde{E}(K, g, \Im, r)$ and $\widetilde{\delta_d}(K, g, \Im ).$ Particular choices of $\widetilde{C_r}$ may have influences on $\widetilde{\delta^i}(K, g, \Im , r)$ and $\widetilde{E}(K, g, \Im, r)$. However, $\widetilde{\delta_d}(K, g, \Im )$ should keep invariant. Therefore, we make the following assumption.

\medskip\noindent
{\bf A Basic Assumption.} {\it When $\delta_d(K, g, \Im )$ does exist for solid $K$ with gravity $g$, small changes of the container boundaries has no effect on $\delta_d(K, g, \Im )$. In other words, for all $\widetilde{C_r}$ satisfying $C_r\subseteq \widetilde{C_r}
\subseteq C_{r+d}$ we have}
$$\widetilde{\delta_d}(K, g, \Im )=\delta_d(K, g, \Im ).$$

It is well-known that, in any given convex domain pack circular disks (different sizes), the density can arbitrarily close to one. Similarly one can show that, in any given convex body vertically pack circular cylinders of shape $C_r$ (different sizes), the density can arbitrarily close to one as well. Thus, when $\Im $ consists only
pouring, shaking and tapping, without rotating (see [14, 15]), the density $\delta_d(K, g, \Im )$ is independent of the shapes of the containers. In other words, we have the following
basic result:

\medskip\noindent
{\bf A Fundamental Theorem.} {\it If $\Im $ is well-defined, the limit $\delta_d(K, g, \Im )$ does exist, and the basic assumption is true, then the density is independent of the shapes of the containers. In other words, let $C$ to be any fixed convex container and replace $C_r$ by $rC$ in the definitions of $\delta^i(K, g, \Im , r)$, $E(K, g, \Im, r)$ and $\delta_d(K, g, \Im )$, the final density $\delta_d(K, g, \Im )$ is invariant.}

\medskip
To achieve the loose random packings (see Scott and Rutgers' experiments), rotation was necessary. Therefore, it seems hard to show that $\delta_l(K,g,\wp )$ is independent of the shapes of the containers.

\medskip
For a fixed convex solid $K$, a set of fixed dense random packing conditions $\Im $,
and a set of fixed loose random packing conditions $\wp$, it is interesting and important to study $\delta_d(K, g, \Im)$
and $\delta_l(K, g, \wp)$ as functions of the specific gravity $g$.

Based on intuitive imagination one may believe that, if $\delta_d(K, g,\Im )$ does exist for some particular specific gravity $g_0$, it should exist for any positive gravity. Further more, simple argument supports the belief that $\delta_d(K, g,\Im )$ should be a non-decreasing stair-like function of the gravity $g$ (constant after sudden change). However, experiments in Dong and Ye \cite{dong93} (see Table 2 and Table 4) contradict the monotonicity. In this setting, some natural questions turn out to be challenging and interesting. For example,

\vspace{0.2cm}
\noindent
{\bf Problem 1.} To determine $\delta_d(B,g,\Im )$, as a function of the specific gravity $g$. In particular, to determine the values of $$\lim_{g\to 0}\delta_d(B, g, \Im ),\quad \lim_{g\to \infty }\delta_d(B, g, \Im ),$$
and
$$\inf_{0<g<\infty }\delta_d(B, g, \Im ), \quad \sup_{0<g<\infty }\delta_d(B, g, \Im ).$$

It is hopeful to determine or estimate the values of $\lim\limits_{g\to 0}\delta_d(B, g, \Im )$ and $\lim\limits_{g\to \infty}\delta_d(B, g, \Im )$ by experiments. For example, apply variant vertical magnetic fields to iron balls can either increase or decrease their gravity. In particular, one can make the gravity as small as possible and therefore obtain estimates on $\lim\limits_{g\to 0}\delta_d(K, g, \Im )$.

Similar questions and considerations can be applied to $\delta_l(K, g, \wp )$ as well, in particular for the ball case.

\medskip
If $\delta_d(K, g, \Im )$ does exist, then
$$\delta_d(\lambda K, g, \Im )=\delta_d(K, g, \Im )$$
holds for any number $\lambda \not= 0$ since $\lambda K$ and $K$ have same shape. Let $\sigma $ denote an affine linear transformation from $\mathbb{E}^3$ to $\mathbb{E}^3$ which keeps the distance ratios invariant. In other words,
$${{\| \sigma ({\bf x}_1), \sigma ({\bf x}_2)\|}\over {\|\sigma ({\bf x}_3), \sigma ({\bf x}_4)\|}}
={{\|{\bf x}_1, {\bf x}_2\|}\over {\| {\bf x}_3, {\bf x}_4\|}}$$
holds for any four points ${\bf x}_1$, ${\bf x}_2$, ${\bf x}_3$ and ${\bf x}_4$, provided ${\bf x}_3\not={\bf x}_4$.
Clearly, $\sigma $ keeps the shape of any geometric object in $\mathbb{E}^3$ unchanged and therefore we have
$$\delta_d(\sigma (K), g, \Im )=\delta_d(K, g, \Im ).$$
It can be easily shown that all such affine linear transformations form a group. For convenience, we denote it by $G$.

Similarly, if $\delta_l(K, g, \wp )$ does exist, for all $\sigma \in G$ we have
$$\delta_l(\sigma (K), g, \wp )=\delta_l(K, g, \wp ).$$

For two convex bodies  $K_1$ and $K_2$ in $\mathbb{E}^3$, let $\gamma $ to be the minimal positive number
such that
$$ K_1\subseteq \sigma (K_2)\subseteq \gamma K_1+{\bf x},\quad \sigma \in G,\quad
{\bf x}\in \mathbb{E}^3,$$
we define
$$\| K_1, K_2||_S=\log \gamma .$$ Clearly, $K_1$ and $K_2$ have same shape if and only if $\gamma =1$ and therefore $\| K_1, K_2\|_S=0$. Furthermore, $K_1$ and $K_2$ have more or less the same shape if and only if  $\| K_1, K_2\|_S$
is small. For convenience, we call $\| K_1, K_2||_S$ the {\it shape distance} between $K_1$ and $K_2$. This metric is similar to the Banach-Mazur metric and Asplund's generalization (see \cite{grub93}), but is neither of them.

Let $B$, $C$ and $T$ denote a ball, a cube and a regular tetrahedron, respectively. One can deduce that
$$\| B, C\|_S={1\over 2}\log 3,\quad \| B, T\|_S=\log 3.$$
On the other hand, as shown by the experiments of Scott \cite{scot60} and Dong and Ye \cite{dong93}, all $\delta_d(B, g, \Im )$, $\delta_d(C, g, \Im )$ and $\delta_d(T, g, \Im )$ do exist, at least for some specific gravities.

\medskip
{\it Does $\delta_d(K, 1, \Im )$ exist for all convex solids $K$}? The answer seems should be \lq\lq no." For example, when $K$ is a cone with a circular disk base of radius 1 cm and height 100 cm, it seems will be hard to achieve real random packings. For this reason, it will be important to characterize all $K$ and $g$ such that $\delta_d(K, g, \Im )$ (or $\delta_l(K, g, \wp )$) does exist and to consider its variants. We end this paper with three problems, which would be important for the completion of this theory and for applications in material science and engineering.

\medskip\noindent
{\bf Problem 2.} With $g$ and $\Im $ are fixed, is $\delta_d(K, g, \Im )$ a continuous function of $K$ with respect to the shape distance?

\medskip\noindent
{\bf Problem 3.} Determine the maximal number $\alpha$ such that $\delta_d(K, g, \Im )$ does exist whenever $\|B, K\|_S\le \alpha .$ In particular, is it true that $\delta_d(K, g, \Im )$ does exist whenever $\| B, K\|_S\le \log 3$?

\medskip\noindent
{\bf Problem 4.} Determine the values of
$$\max_{K, g}\delta_d(K, g, \Im )\quad {\rm and} \quad \min_{K, g}\delta_d(K, g, \Im ),$$
where the maximum and the minimum are over all pairs such that $\delta_d(K, g, \Im )$ does exist.

\vspace{0.5cm}\noindent
{\bf Acknowledgements.} This work is supported by 973 Programs 2013CB834201 and 2011CB302401, the National Natural Science Foundation of China (No. 11071003), and the Chang Jiang Scholars Program of China.

\bibliographystyle{amsplain}

\end{document}